\documentclass[12pt]{article}

\usepackage{amsmath,amssymb,amsfonts,theorem,makeidx,latexsym,epsfig,,subfigure}

\newtheorem{defn}{Definition}[section]

\newtheorem{lemma}[defn]{Lemma}

{\theorembodyfont{\rmfamily}

\newtheorem{ex}[defn]{Example}}

\newtheorem{thm}[defn]{Theorem}

\newtheorem{prop}[defn]{Proposition}

\newtheorem{cor}[defn]{Corollary}

\newtheorem{rem}[defn]{Remark}

\numberwithin{equation}{section}

\newcommand{\h}{{\cal H}}

\newcommand{\ltr}{ L^2(\mathbb R) }

\newcommand{\ltn}{{\ell}^2(\mathbb N)}

\newcommand{\ltz}{{\ell}^2(\mathbb Z)}

\newcommand{\si}{S^{-1}}

\newcommand{\mn}{\mathbb N}

\newcommand{\mr}{\mathbb R}

\newcommand{\mz}{\mathbb Z}

\newcommand{\mc}{\mathbb C}

\newcommand{\mts}{ \{E_{mb}T_{na}g \}_{m,n \in \mz}}

\def\bp{{\noindent\bf Proof. \ }}

\def\ep{\hfill$\square$\par\bigskip}

\def\bqs{\begin{equation}}

\def\eqs{\tag*{$\square$}\end{equation}\par\bigskip}

\def\la{\langle}

\def\ra{\rangle}

\def\ftk{\{f_k\}_{k=1}^\infty}

\def\ctk{\{c_k\}_{k=1}^\infty}

\def\gtk{\{g_k \}_{k=1}^\infty}

\def\htk{\{h_k\}_{k=1}^\infty}

\def\etk{\{e_k\}_{k=1}^\infty}

\def\suk{\sum_{k=1}^\infty}

\def\nl{\left|\left|}

\def\nr{\right|\right|}

\def\span{\overline{\text{span}}}

\def\Span{\text{span}}

\def\supp{\text{supp}}

\def\vn{\vspace{.1in}\noindent}

\def\bop{\begin{op}\rm}

\def\eop{\end{op}}

\def\bee{\begin{eqnarray}}

\def\ene{\end{eqnarray}}

\def\bes{\begin{eqnarray*}}

\def\ens{\end{eqnarray*}}

\def\bei{\begin{itemize}}

\def\eni{\end{itemize}}

\def\bt{\begin{thm}}

\def\et{\end{thm}}

\def\bc{\begin{cor}}

\def\ec{\end{cor}}

\def\bpr{\begin{prop}}

\def\epr{\end{prop}}

\def\bl{\begin{lemma}}

\def\el{\end{lemma}}

\def\bd{\begin{defn}}

\def\ed{\end{defn}}

\def\bex{\begin{ex}}

\def\enx{\end{ex}}

\def\bfi{\begin{fig}}

\def\efi{\end{fig}}

\def\Tnn{\{T^n\varphi\}_{n=0}^\infty}

\title{Approximate frame representations via   iterated operator systems}

\date{\today}

\author{Ole Christensen, Marzieh Hasannasab }

\begin{document}

\maketitle

\begin{abstract}
	It is known that it is a very restrictive condition for a frame $\ftk$ to
have a representation $ \{T^n \varphi\}_{n=0}^\infty$ as the orbit of
a bounded operator $T$  under a single generator $\varphi \in \h.$  In this paper
we prove that, on the other hand, any frame can be approximated arbitrarily well by a suborbit
$\{T^{\alpha(k)} \varphi\}_{k=1}^\infty$ of a bounded operator $T$. An important
new aspect is that for certain important classes of frames, e.g., frames
consisting of finitely supported vectors in $\ltn,$ we can be
completely explicit about possible choices  of the operator $T$ and
the powers $\alpha(k), \, k\in \mn.$ A similar approach carried out
in $\ltr$ leads to an approximation of a frame using suborbits of two
bounded operators. The results are illustrated with an
application to  Gabor
frames generated by a compactly supported function. The paper is concluded with an appendix which collects general
results about frame representations using multiple orbits of bounded
operators.
\end{abstract}

\section{Introduction} \label{82309a}

In the recent and very active research area dynamical sampling
\cite{A1,A2,A3,CMPP,olemmaarzieh-2,olemmaarzieh-3}, one of the key issues is to consider frames of the form
$ \{T^n \varphi\}_{n=0}^\infty$, where $T$ is a bounded linear
operator on a given Hilbert space $\h$ and $\varphi \in \h.$
Unfortunately, the class of explicitly known frames with
such  a representation  is very sparse: except for the
Riesz bases, it only contains the Carleson frames introduced in \cite{A1}. Various alternative operator representations can be found in the literature, but each of these come with its own
constraints and limitations. The purpose of this paper is to
show that all of these limitations disappear if we consider
approximate frame representations using suborbits of a bounded
operator. Indeed, we prove that any frame $\ftk$ can be approximated arbitrarily well in a perturbation theoretic sense by a frame $\{T^{\alpha(k)} \varphi\}_{k=1}^\infty,$ which
shares key features with the given frame. Appropriate choices
of the operator $T$ already exist in the literature, e.g., the
hypercyclic operators;  the most important new aspects in the current paper is that we for some important classes of frames, e.g., frames consisting of finitely supported
vectors in $\ltn,$
also can specify appropriate choices of the powers $\alpha(k), \, k\in \mn.$
A generalization to vectors that are not finitely supported but with coordinates
``decaying sufficiently fast" is presented as well.
We also treat the case of frames consisting of compactly supported functions in
$\ltr,$ which for technical reasons leads to an approximation using
suborbits of two bounded operators rather than one. The results are illustrated with an
application to  Gabor
frames generated by a compactly supported function. Since the idea of using
multiple suborbits is not treated detailed in the existing literature,
the paper is concluded with an appendix, which discusses various frame representations
using multiple suborbits and the constraints under which the results are applicable.

The paper is organized as follows. In the rest of this section
we set the stage by discussing the known results and limitations
concerning represesentations of frames as an orbit of a bounded operator. The results about approximate frame representations are
in Section \ref{82309e}.

The following result collects some of the results from the
literature about various operator representations of frames. In the entire paper
we let $\h$ denote an infinite-dimensional separable Hilbert space. Recall
that a frame $\ftk$ is said to have {\it infinite excess} if there
exists an infinite subfamily $\{f_k\}_{k\in I}$ such that
$\{f_k\}_{k\in \mn \setminus I}$ is a frame; it has {\it finite excess}
if it consists of a Riesz basis and a finite number of vectors.

\bt \label{80923f}   Let $\ftk$ be a frame for $\h.$ Then the following hold.
\bei
\item[(i)]  {\bf \cite{olemmaarzieh,CHP}} Assume that $\ftk$ is overcomplete and has a representation
$\ftk= \Tnn$ for a bounded operator $T:\h \to \h.$  Then
$\ftk$ is linearly independent,
has infinite excess, and $f_k\to 0$ as $k\to \infty.$
\item[(ii)] {\bf \cite{kitai}} Assume that $\{f_k\}_{k=1}^\infty$ is  linearly independent. Then, given any   $a>1$,
there exists  a bounded linear operator $T: \h \to \h$ with $\| T \| = a$ and a sequence $\{\alpha(k)\}_{k=1}^\infty\subset\mn_0$ such that $\{f_k\}_{k=1}^\infty=\{T^{\alpha(k)}f_1\}_{k=1}^\infty$.
\eni
\et

It is worth noticing that the conditions in
Theorem \ref{80923f}(i) are only necessary but not sufficient for the availability of a representation of a frame
as an orbit of a bounded operator. Various characterizations of frames with such a representation can be found in
\cite{olemmaarzieh-3} and \cite{CHP}. A generalization of Theorem \ref{80923f}(ii) to Banach spaces appeared in \cite{Grivaux}.

The proof of Theorem \ref{80923f}(ii) is given in terms of a procedure and is not easy to implement.
Our results in Section \ref{82309e} are indeed inspired by
Theorem \ref{80923f}(ii) and the desire to be more
explicit about the operator $T$ and the sequence $\{\alpha(k)\}_{k=1}^\infty\subset\mn_0$.
Note that the result in Theorem \ref{80923f}(ii) actually does not
need the frame property. On the other hand, in the
current context the linear independency is indispensable:

\bpr\label{07112017a} Assume that $\{T^{\alpha(k)}\varphi\}_{k=1}^\infty$ spans an infinite-dimensional space
and that $\alpha(k) \neq \alpha(\ell)$ for $k\neq \ell.$ Then $\{T^{\alpha(k)}\varphi\}_{k=1}^\infty$ is linearly independent.
\epr

\bp Without loss of generality we can assume that the vectors  $\{T^{\alpha(k)}\varphi\}_{k=1}^\infty$ are ordered such that the sequence
$\{\alpha(k)\}_{k=1}^\infty$ is increasing. Now,
if $T^{\alpha (N)} \varphi \in \Span \{T^{\alpha(k)}\varphi\}_{k=1}^{N-1}$
for some $N \ge 2,$ then the vector space
$$V:= \span \{\varphi, T\varphi, \dots, T^{\alpha(N)-1}\varphi\}$$ is invariant under the operator $T$ and hence
$\{T^{\alpha(k)}\varphi\}_{k=1}^\infty$ only spans a finite-dimensional space.\ep

Note that
Theorem \ref{80923f}(i) implies that the standard frames in harmonic analysis, e.g.,
overcomplete frames of translates, Gabor
frames and wavelet frames, do not have representations as a (full) orbit of a bounded
operator. However, by Theorem \ref{80923f}(ii) a representation as a suborbit of
a bounded operator is available under the condition of linear independency, a condition
that is satisfied for all regular frames of translates, Gabor frames along
lattices \cite{Linn}, and a large class of wavelet frames \cite{BoSp2}.

\section{Approximate
frame representations}\label{82309e}
The existing literature on dynamical sampling deals with exact frame representations.
The price to pay for this is that although the class of frames that can be
represented as an orbit of a bounded operator is quite large \cite{CHP}, only
few explicitly given frames of that form are known \cite{A1,olemmaarzieh}. In this section we
show that all constraints disappear if we consider
approximate frame representations. Indeed, we prove that any
frame can be approximated up to any desired precision
using a frame of the form $\{T^{\alpha(k)}\varphi\}_{k=1}^\infty$ for certain
choices of the operator $T,$ the generator $\varphi\in \h,$ and the powers $\alpha(k)\in \mn.$ The ``approximation quality"  is measured in terms  of classical conditions from
perturbation theory for frames, which are known to preserve key properties of the frame, e.g., the excess.  For the frames
of particular relevance in applications, e.g.,
frames consisting of  finitely supported vectors in $\ltn$ and compactly supported
functions in $\ltr,$ we can specify exact choices of
the powers $\alpha(k)$ for certain choices of the operator
$T$ and the generator $\varphi.$ The general perturbation theoretic approach is outlined in Section \ref{191018a},
and the applications to frames in $\ltn$ and $\ltr$ are
given in Sections \ref{191018b}--\ref{191018c}.

\subsection{Perturbation theory and hypercyclic operators} \label{191018a}

For later convenience we formulate the definition below for
a frame with an arbitrary countable index set $I.$

\bd \label{191001a} Let $\{f_k\}_{k\in I}$ be a frame for $\h.$ Given any $\epsilon >0,$
a sequence $\{\widetilde{f_k}\}_{k\in I} \subset \h$ is called an $\epsilon$-approximation of
$ \{f_k\}_{k\in I}$ if
\bee \label{151018a}  \nl \sum c_k (f_k - \widetilde{f_k})\nr^2
\le \epsilon \, \sum |c_k|^2\ene
for all finite sequences $\{c_k\}.$\ed

For sufficiently small values of $\epsilon,$ an $\epsilon$-approximation $\{\widetilde{f_k}\}_{k\in I}$ of a frame $\{f_k\}_{k\in I}$
is itself a frame and shares several key properties
of the frame, e.g., its excess. Furthermore, the synthesis
operator and the frame operator for $\{f_k\}_{k\in I}$ are
approximated ``well" by the corresponding operators for
$\{\widetilde{f_k}\}_{k\in I}$:

\bt\label{181217ens}  Consider a frame $\{f_k\}_{k\in I}$
for $\h$  with frame bounds $A,B$, and assume that
$\{\widetilde{f_k}\}_{k\in I} \subset \h$ is an $\epsilon$-approximation of
$ \{f_k\}_{k\in I}$ for some $\epsilon \in ]0, A[.$
Then the following hold:

\bei\item[(i)]  $ \{\widetilde{f_k}\}_{k\in I} $ is a frame with bounds $A(1-\sqrt{ \frac{\epsilon}{A}  })^2$ and $B(1+\sqrt{ \frac{\epsilon}{B}  })^2,$ with the same excess as
$ \{f_k\}_{k\in I}.$
\item[(ii)] Denoting the synthesis operators and frame operators of $\{f_k\}_{k\in I}$ and
    $\{\widetilde{f_k}\}_{k\in I}$  by $U, \widetilde{U}$, respectively, $S,\widetilde{S}$, we have \[ \| U- \widetilde{U} \| \leq \sqrt{\epsilon},\quad
  \| S - \widetilde{S} \|\leq \sqrt{\epsilon B}  \left(2+\sqrt{ \frac{\epsilon}{B}  }  \right),\]
  and \[
    \| S^{-1} - \widetilde{S}^{-1} \|\leq  \frac{ \sqrt{\epsilon B} (2+\sqrt{ \frac{\epsilon}{B}}) }{A^2(1-\sqrt{ \frac{\epsilon}{A}  })^2}. \]
\eni
\et
\bp $(i)$ follows from Corollary 22.1.5 in \cite{CB}.
For the proof of $(ii)$, letting $\ctk\in\ltn,$
\bes \| U \ctk - \widetilde{U} \ctk \|^2 &=& \| \suk c_k (f_k - \widetilde{f_k}) \|^2 \\
&\leq & \epsilon \suk | c_k |^2. \ens
So $\| U- \widetilde{U} \| \leq \sqrt{\epsilon}. $ By the definition of frame operator, we have
\bes S- \widetilde{S} = UU^* - \widetilde{U} \widetilde{U}^*= (U - \widetilde{U} )U^* + \widetilde{U}  (U^* - \widetilde{U}^*). \ens It is well known that if $B$ is any
upper frame bound for a frame with synthesis operator $U,$
then $\| U \| \leq \sqrt{B}.$ Thus
it follows from (i) that  $\| \widetilde{U} \| \leq \sqrt{B(1+\sqrt{\frac{\epsilon}{B}})^2}$. Therefore
\[ \| S - \widetilde{S} \| \leq \sqrt{\epsilon} \sqrt{B} + \sqrt{\epsilon} \sqrt{B
	\left(1+\sqrt{ \frac{\epsilon}{B}  }  \right)^2  }=\sqrt{\epsilon B}   \left(2+\sqrt{ \frac{\epsilon}{B}  }  \right).\]
Also, if $A$ is any lower frame bound
for a frame with frame operator $S,$ then $||\si|| \le A^{-1}.$
For the inverse of the frame operators we therefore obtain that
 \bes \| S^{-1} - \widetilde{S}^{-1} \| &=& \| S^{-1}(S-\widetilde{S}) \widetilde{S}^{-1} \| \leq  \frac{1 }{A}  \sqrt{\epsilon B} \left(2+\sqrt{ \frac{\epsilon}{B}}\right) \frac{ 1 }{A(1-\sqrt{ \frac{\epsilon}{A}  })^2},
 \ens  as claimed. \ep

For the applications in Sections \ref{191018b}-\ref{191018c}
it is convenient to apply the following sufficient condition
for $\{\widetilde{f_k}\}_{k\in I}$ being an $\epsilon$-approximation
of $\{f_k\}_{k\in I}$:

\bpr \label{231018e} Let $\{f_k\}_{k\in I}$ be a frame for $\h$
with lower frame bound $A$, and assume
that for some $\epsilon \in [0, A[$ the sequence $\{\widetilde{f_k}\}_{k\in I} \subset \h$
satisfies that
\bee \label{151018b} \sum_{k\in I} ||f_k- \widetilde{f_k}||^2
\le \epsilon.\ene
Then $\{\widetilde{f_k}\}_{k\in I}$ is an $\epsilon$-approximation
of $\{f_k\}_{k\in I}$ and the conclusions in Theorem \ref{181217ens} hold. \epr

The applications in the rest of the section will
be based on the following special case of
Proposition \ref{231018e}:

\bc\label{181217en}  Let $\ftk$ be a frame
for $\h$  with lower frame bound $A$: Let  $\varphi \in \h$
and consider a
bounded operator $T: \h \to \h$. Assume
that  for a given $\epsilon \in ]0, A[,$ and for any $k\in \mn$ there exists a nonnegative integer
$\alpha(k)\in\mn_0$ such that
\bee \label{011018p}  \| f_k -  T^{\alpha(k)} \varphi \|^2 \leq\frac{\epsilon}{2^k}.\ene
Then $\{ T^{\alpha(k)} \varphi \}_{k=1}^\infty$ is an $\epsilon$-approximation
of $\ftk$ and the conclusions in Theorem \ref{181217ens} hold.
\ec

Corollary \ref{181217en} is inspired by the theory for hypercyclic
operators.
Recall that a bounded operator $T: \h \to \h$ is
{\it hypercyclic} if there exists a vector $\varphi\in \h$ such that the orbit $\{T^n\varphi\}_{n=0}^\infty$
is dense in $\h;$ any such vector $\varphi$ is called a
hypercyclic vector. We refer to, e.g., \cite{erdmann,erdmannsurvey} for more
information on hypercyclic operators.

\bex \label{171218a} Let $T:\h \to \h$ be a hypercyclic operator with hypercyclic vector $\varphi.$  Then, for any frame
$\ftk$ for $\h$ and any given $\epsilon >0$ there exists
nonnegative integers $\alpha(k), \, k\in \mn,$ such that
\eqref{011018p} holds. Thus, for sufficiently small
values of $\epsilon,$ we obtain a frame $\{T^{\alpha(k)}\varphi\}_{k=1}^\infty$ with the same
excess as the frame $\ftk$ and approximating $\ftk$
in the sense of Definition \ref{191001a}.
\ep \enx

It is worth emphasizing that in Example \ref{171218a}
the obtained frames $\{T^{\alpha(k)}\varphi\}_{k=1}^\infty$
are only suborbits of the orbit for the operator $T.$
Indeed, the full orbit of a hypercyclic operator
can never be a frame, regardless of the choice of
the initial vector:

\bpr Assume that $T: \h\to \h$ is hypercyclic. Then $\{T^n y\}_{n=0}^\infty$
and $\{(T^*)^n y\}_{n=0}^\infty$ are not  frames for
any choice of $y \in \h.$ \epr

\bp Let $\varphi\in \h$ be a hypercyclic vector and consider any $y\in \h \setminus \{0\}.$ Then the sequence $\{\la y, T^n \varphi\ra\}_{n=0}^\infty =
\{ \la (T^*)^ny, \varphi\ra\}_{n=1}^\infty$ is unbounded, which implies that
$\{||(T^*)^n y||\}_{n=1}^\infty$ is unbounded. Thus, by Theorem 7 in \cite{A3} the
sequence $\{T^n y\}_{n=0}^\infty$ is not a frame.
Also, considering again any $y \neq 0,$
\bes \sum_{n=0}^\infty | \la \varphi, (T^*)^ny\ra|^2=
\sum_{n=0}^\infty | \la T^n\varphi, y\ra|^2=\infty,\ens
which implies that $\{(T^*)^n y\}_{n=0}^\infty$ is not a frame. \ep

The first example of a hypercyclic operator was
given by Rolewicz in \cite{rolewicz}, who showed that
for any $\lambda>1$ the scaled left-shift operator
$T(x_1, x_2, x_3, \dots):= \lambda (x_2, x_3, \dots)$ is
hypercyclic on $\ltn.$ More recently a hypercyclic operator
has been constructed on $\ltr$:

\bex \label{231018a} Let $w: \mr \to \mr$ be a continuous decreasing function
such that
\bes w(x)= \begin{cases} 2 \, &\mbox{if} \, \, x\le 0 \\
1/2  \, \, &\mbox{if} \, x\ge 1.  \end{cases} \ens
Then the multiplication operator $T: \ltr \to \ltr,
Tf(x):= w(x)f(x)$ is hypercyclic, \cite{Chen}. In particular, reindexing
any frame for $\ltr$ as $\ftk,$ we can approximate the
frame arbitrarily well in the sense of Definition \ref{191001a}
via an appropriately chosen suborbit $\{T^{\alpha(k)}\varphi\}_{k=1}^\infty.$\ep \enx

Example \ref{231018a} points directly to one of the key issues
of the current paper: indeed, in general we do not have direct
access to the  numbers  $\alpha(k), \, k\in \mn,$ such that
\eqref{011018p} holds, not even for
hypercyclic operators. The purpose of the next section is
to show that under very natural conditions on frames in
$\ltn$ and $\ltr$ we can be very explicit about the choices
of appropriate scalars $\alpha(k)$ for certain choices
of the operator $T.$

\subsection{Approximate frame representations in $\ltn.$}
\label{191018b}

In this section we will consider a frame $\ftk$ for $\ltn$
consisting of finitely supported vectors. In the main result we consider a particular
operator $T$ and show how to choose a vector $\varphi\in \ltn$ and powers $\alpha(k)$ such that \eqref{011018p} is
satisfied. In order to avoid a too cumbersome formulation of
the results we will apply the following standing assumptions
throughout the section:

\vn{\bf General setup:} In the entire section, let $\ftk$ denote a frame for $\ltn$
with frame bounds $A,B.$ Fix some $\lambda >1$, and consider the scaled
left/right-shift operators on $\ltn,$ given by
\bee \label{231018d} T(x_1, x_2, \dots)=\lambda(x_2, x_3,\dots), \, \, \,
U(x_1, x_2, \dots)=\lambda^{-1}(0,x_1, x_2,\dots).\ene
For appropriately chosen nonnegative integers $\alpha(k), \, k\in \mn,$
let
\bee \label{011018bn}
\varphi:= \sum_{n=1}^\infty U^{\alpha(n)} f_n.\ene

For the case of finitely supported vectors $\ftk$
the following result specifies how to choose the powers
$\alpha(k)$ such that the vector $\varphi$ in \eqref{011018bn}
is well-defined and \eqref{011018p} holds:

\bt \label{011018an} Under the conditions in the general setup, assume that $\ftk$  consists of finitely supported vectors
in $\ltn;$ for
$k\in \mn,$ let $m(k)$ denote the largest index for a nonzero
coordinate in $f_k.$ Let $\{\alpha(k)\}_{k=1}^\infty$ be a
strictly increasing sequence
of nonnegative integers such that $\alpha(1)=0$ and
$\alpha(k+1)-\alpha(k)\ge m(k)$ for all $k\in \mn.$
Then,   for any  $k\in \mn,$
\bee \label{011018cn} || f_k- T^{\alpha(k)}\varphi||^2 \le  \frac{B\lambda^2}{\lambda^2-1} \lambda^{-2[\alpha(k+1)-\alpha(k)]}.\ene
In particular, choosing the nonnegative integers $\alpha(k)$ such that $\alpha(1)=0$ and
for a given $\epsilon \in ]0, A[,$
\bee \label{011018k}
\alpha(k+1)-\alpha(k) \ge \max\left(m(k), \frac{k \,\ln(2)+\ln \left(\frac{B}{\epsilon}\right) + \ln \left( \frac{\lambda^2}{\lambda^2-1}\right)  }{2\ln (\lambda)}\right), \ene
the condition \eqref{011018p} is satisfied, i.e., the conclusions in Theorem \ref{181217ens} hold.

\et

\bp We have $\supp\, f_n \subseteq \{1, \dots, m(n)\}$ and
hence $$\supp\, U^{\alpha(n)}f_n \subseteq \{1+ \alpha(n), \dots, m(n)+\alpha(n)\};$$
thus the assumption $\alpha(k+1)-\alpha(k)\ge m(k), \, \forall k\in \mn,$ implies
that the vectors $U^{\alpha(n)}f_n, \, n\in \mn,$ are
perpendicular.
Using that $\{\alpha(k)\}_{k=1}^\infty$ is assumed to be strictly increasing and that the operator $U$ is a multiple of an isometry, it follows that
\bee \label{031218a} \sum_{n=1}^\infty ||U^{\alpha(n)} f_n||^2=
\sum_{n=1}^\infty \lambda^{-2\alpha(n)} ||f_n||^2
\le B\sum_{\ell=0}^\infty \lambda^{-2\ell} = \frac{ \lambda^2 \, B}{\lambda^2 -1} < \infty.\ene
Thus the infinite sequence defining the vector $\varphi$
is convergent. Now, fix any $k\in \mn.$ Then
$T^{\alpha (k)} U^{\alpha(k)}f_k=f_k$. Also, for any $n <k,$ we have
$ T^{\alpha (k)} U^{\alpha(n)}f_n=  T^{\alpha(k)- \alpha(n)}f_n;$ since
\bes \alpha(k)-\alpha(n) \ge \alpha(n+1)-\alpha(n) \ge m(n),\ens
it follows that $T^{\alpha (k)} U^{\alpha(n)}f_n=0$ for $n<k.$
Thus
\bee \label{031218b} ||f_k - T^{\alpha(k)}\varphi ||^2 & = &
\nl \sum_{n=k+1}^\infty T^{\alpha(k)} U^{\alpha(n)} f_n \nr^2
\\ \notag & = & \nl \sum_{n=k+1}^\infty  U^{\alpha(n)-\alpha(k)} f_n \nr^2
\\ \notag & = &   \sum_{n=k+1}^\infty \lambda^{-2\, [\alpha(n)-\alpha(k)]} ||f_n ||^2 \\ \notag & \le &
B \lambda^{-2\, [\alpha(k+1)-\alpha(k)]} \sum_{\ell=0}^\infty
\lambda^{-2\ell} \\ \notag & = & \frac{B\lambda^2}{\lambda^2-1} \lambda^{-2[\alpha(k+1)-\alpha(k)]}.\ene
This proves \eqref{011018cn}. In order to satisfy
the condition \eqref{011018p} it is now enough
to put the additional constraint on $\alpha(k)$
that
\bes \frac{B\lambda^2}{\lambda^2-1} \lambda^{-2[\alpha(k+1)-\alpha(k)]} \le \epsilon 2^{-k};\ens this leads to the sufficient condition stated in
\eqref{011018k}. \ep

The inequality \eqref{011018k} directly express how
the powers $\alpha(k)$ are influenced
by the parameter $\lambda,$ the support lengths $m(k)$, and
the desired precision level determined by the parameter $\epsilon.$ In the next result we derive an exact
expression for suitable powers $\alpha(k)$ for
the special case $\lambda=\sqrt{2}.$

\bc \label{180226a} In the setup of Theorem \ref{011018an}, let $\lambda=\sqrt{2},$ take an upper frame bound of the
form $B=2^N$ for some $N\in \mn$ and a tolerance
$\epsilon =2^{-j}$ for some $j\in \mn.$ Then the following hold.

\bei
\item[(i)] Without any restriction on the support sizes
$m(k)$ of the vectors $f_k,$ the condition \eqref{011018k}
is satisfied if $\alpha(1)=0$ and
\bee \label{061018a} \alpha(k)= (k-1)\left[ N+j + 1+\frac{k}{2}\right] + \sum_{\ell=1}^{k-1} m(\ell), \, \, k\in \mn \setminus\{1\}. \ene\item[(ii)] If
$m(k) \le N+j + 1+k$ for all $k\in \mn,$ the condition \eqref{011018k}
is satisfied if
\bee \label{061018b} \alpha(k)= (k-1)\left[ N+j + 1+\frac{k}{2}\right], \, \, k\in \mn. \ene \eni

\ec

\bp For the given choice of the parameters $\lambda, B, \epsilon,$ a direct calculation gives that
\bes \frac{k \,\ln(2)+\ln \left(\frac{B}{\epsilon}\right) + \ln \left( \frac{\lambda^2}{\lambda^2-1}\right)  }{2\ln (\lambda)}= k+N+j+1;\ens thus, without any assumption
on the numbers $m(k),$ the condition \eqref{011018k}
is satisfied if
\bee \label{061018c}
\alpha(k+1)-\alpha(k) = k+N+j+1+ m(k).\ene Using that
$\alpha(1)=0$ this yields the formula \eqref{061018a}.
In case $m(k) \le N+j + 1+k$ for all $k\in \mn,$ we
can discard the numbers $m(k)$ in \eqref{061018c}, which then
leads to the formula stated in \eqref{061018b}. \ep

Under certain decay conditions on the coordinates
of the vectors $f_k,$ we can remove the assumption of
finite support in Theorem \ref{011018an}. Specifically,
letting $\etk$ denote the canonical orthonormal basis
for $\ltn,$ the $j$th coordinate of the vector
$f_k$ is $\la e_j,f_k\ra;$ we will assume
that that there exist constants  $C,\beta >0$ such that
\bee \label{041018a} | \la e_j, f_k\ra| \le C e^{-\beta \, |j-k|}, \, \forall j,k \in \mn.\ene
Note that in the frame literature this is phrased by saying that the
frame $\ftk$ is
{\it $\beta$-exponentially localized} with respect to the orthonormal basis
$\etk$; see \cite{G5}.

\bt \label{011018am} Under the conditions in the
general setup, assume that the frame
$\ftk$ is $\beta$-exponentially
localized with respect to the canonical orthonormal
basis $\etk$ for $\ltn,$  as in \eqref{041018a}.
Consider an increasing sequence $\{\alpha(k)\}_{k=1}^\infty$
of nonnegative integers such that $\alpha(1)=0$ and
$\alpha(k)\ge \alpha(k-1)+k-2$ for all $k\ge 2.$
Then
\bee \label{041018b} ||f_1- T^{\alpha(1)}\varphi|| \le
\frac{\lambda}{\lambda-1} B^{1/2} \lambda^{-\alpha(2)},\ene
and, for any  $k\in \mn \setminus\{1\},$
\bee \notag   || f_k- T^{\alpha(k)}\varphi||  & \le &
\frac{\lambda}{\lambda-1} B^{1/2} \lambda^{-[\alpha(k+1)-\alpha(k)]} \\ \label{011018cm}
& \ & + (\lambda e^{-\beta})^{\alpha(k)}\left(
\frac{Ce^{-\beta}}{\sqrt{1-e^{-2\,\beta}}}  \sum_{n=1}^{k-1} (\lambda e^{-\beta})^{-\alpha(n)}e^{\beta n}\right). \, \hspace{1cm} \,
\ene

In particular, assuming that $\ln(\lambda)< \beta$ and
fixing any $\epsilon \in ]0,A[$ if we choose the nonnegative integers $\alpha(k)$
recursively  such that $\alpha(1)=0$ and for
$k\ge 1,$

\bee \notag \alpha(k+1) \ge\max\left( \alpha(k)+k-1,   \alpha(k)+  \frac{(k/2+1) \ln 2   + \ln \frac{\lambda}{\lambda-1} + \ln \sqrt{\frac{B}{\epsilon}}}{\ln \lambda}, \right. \\ \label{011018b}
\left. \frac{(k/2+3/2)  \ln 2 + \ln \left(\sum_{n=1}^{k} (\lambda e^{-\beta})^{-\alpha(n)}e^{\beta n}\right) + \ln \left(
 \frac{Ce^{-\beta}}{\sqrt{1-e^{-2\,\beta}}}
    \right) -\ln \sqrt{\epsilon} }{\beta - \ln \lambda}\right), \ene
the condition \eqref{011018p} is satisfied, i.e., the conclusions in Theorem \ref{181217ens} hold.

\et

\bp We first notice that the infinite series defining $\varphi$
in \eqref{011018bn} is convergent; indeed, since the sequence
$\{\alpha(k)\}_{k=1}^\infty$ is assumed to be strictly increasing and the operator $U$ is a multiple of an isometry,
\bes \sum_{n=1}^\infty ||U^{\alpha(n)} f_n||=
\sum_{n=1}^\infty \lambda^{-\alpha(n)} ||f_n||
\le \sqrt{B}\sum_{\ell=0}^\infty \lambda^{-\ell} = \frac{ \lambda \, \sqrt{B}}{\lambda -1} < \infty.\ens

Now, using that $\alpha(1)=0,$ the same calculation yields that
\bes ||f_1-T^{\alpha(1)}\varphi|| = \nl \sum_{n=2}^\infty U^{\alpha(n)}f_n \nr \le \sum_{n=2}^\infty \lambda^{-\alpha(n)} ||f_n||\le \frac{\lambda}{\lambda-1} B^{1/2} \lambda^{-\alpha(2)},\ens as claimed in \eqref{041018b}.
Considering now any $k\ge 2$ and using that
$f_k= T^{\alpha(k)}U^{\alpha(k)}f_k,$
\bee \notag \label{041018c}  || f_k- T^{\alpha(k)}\varphi||
& = & \nl \sum_{n=1}^{k-1} T^{\alpha(k)}U^{\alpha(n)}f_n  + \sum_{n=k+1}^\infty T^{\alpha(k)}U^{\alpha(n)}f_n\nr \\ & \le &
\sum_{n=1}^{k-1} \nl T^{\alpha(k)-\alpha(n)}f_n\nr +
\sum_{n=k+1}^\infty \nl U^{\alpha(n)-\alpha(k)}f_n\nr.
\ene We now consider the two terms in \eqref{041018c} separately. For the first first term, note that since
$\alpha(k) > \alpha(n)$ for $n=1, \dots, k-1,$ the operator
$T^{\alpha(k)-\alpha(n)}$ removes the first $\alpha(k)-\alpha(n)$ coordinates of the vector the operator
is acting on, and multiplies the resulting vector with
$\lambda^{\alpha(k)-\alpha(n)}.$ Thus
\bee \notag  \nl T^{\alpha(k)-\alpha(n)}f_n\nr
& = &   \lambda^{\alpha(k)-\alpha(n)}
\left( \sum_{j=  \alpha(k)-\alpha(n)+1}^\infty | \la e_j, f_n\ra|^2\right)^{1/2} \\ & \le &
  \lambda^{\alpha(k)-\alpha(n)}
 \left( \sum_{j=  \alpha(k)-\alpha(n)+1}^\infty C^2\, e^{-2\beta \, |j-n|}\right)^{1/2}. \label{041018d}
 \ene The assumption $\alpha(k)\ge \alpha(k-1)+k-2$ implies
 that $j-n \ge 0$ for all
 $n=1, \dots, k-1$ and all $j\ge \alpha(k)-\alpha(n)+1;$
 thus,  \eqref{041018d} yields that
 \bee \notag \nl T^{\alpha(k)-\alpha(n)}f_n\nr & \le &
 \lambda^{\alpha(k)-\alpha(n)}
 \left( \sum_{j= 0}^\infty C^2\, e^{-2\beta \, (\alpha(k)-\alpha(n)+1+j-n)}\right)^{1/2} \\ \label{041018f}
 & = &  (\lambda e^{-\beta})^{\alpha(k)}\frac{Ce^{-\beta}}{\sqrt{1-e^{-2\,\beta}}} (\lambda e^{-\beta})^{-\alpha(n)}e^{\beta n}
 \ene  Considering now the second term in \eqref{041018c},
 \bee \notag \sum_{n=k+1}^\infty \nl U^{\alpha(n)-\alpha(k)}f_n\nr
 & = &   \sum_{n=k+1}^\infty \lambda^{-(\alpha(n)-\alpha(k))} \nl f_n   \nr \\ & \le & \notag
B^{1/2} \lambda^{-(\alpha(k+1)-\alpha(k))} \sum_{\ell=0}^\infty \lambda^{-\ell} \\ & \le & \frac{\lambda}{\lambda-1} B^{1/2} \lambda^{-[\alpha(k+1)-\alpha(k)]}.  \label{041018g}
 \ene Now \eqref{011018cm} follows immediately from
\eqref{041018c} using the expressions in \eqref{041018f}
and \eqref{041018g}.

Now, using that $\alpha(1)=0,$ the inequalities \eqref{041018b} and \eqref{011018cm} show that \eqref{011018p} holds if
we choose $\alpha(k), k \ge 2,$ such that
\bee \label{221018a}
\frac{\lambda}{\lambda-1} B^{1/2} \lambda^{-[\alpha(k+1)-\alpha(k)]} \le \frac12 \sqrt{\epsilon 2^{-k}}\ene and
\bee \label{221018b}
(\lambda e^{-\beta})^{\alpha(k)}\left(
\frac{Ce^{-\beta}}{\sqrt{1-e^{-2\,\beta}}}  \sum_{n=1}^{k-1} (\lambda e^{-\beta})^{-\alpha(n)}e^{\beta n}\right)
 \le \frac12 \sqrt{\epsilon 2^{-k}}. \ene
 Direct calculations (which we skip) yield that \eqref{221018a}
 is satisfied if
 \bes \alpha(k+1)-\alpha(k) \ge
 \frac{(k/2+1) \ln 2   + \ln \frac{\lambda}{\lambda-1} + \ln \sqrt{\frac{B}{\epsilon}}}{\ln \lambda}
 \ens
 and that \eqref{221018b} holds if
 \bes \alpha(k) \ge \frac{(k/2+1)  \ln 2 + \ln \left(\sum_{n=1}^{k-1} (\lambda e^{-\beta})^{-\alpha(n)}e^{\beta n}\right) + \ln \left(
 \frac{Ce^{-\beta}}{\sqrt{1-e^{-2\,\beta}}}
    \right) -\ln \sqrt{\epsilon} }{\beta - \ln \lambda}.\ens
This completes the proof of \eqref{011018b}.\ep

\begin{rem} \label{231018c} {\rm The approach taken in this section is not restricted to the Hilbert space $\ltn.$ In fact, taking any separable Hilbert space $\h$ and an orthonormal basis $\etk$, all the
stated results can be formulated in the setting of
the Hilbert space $\h$ using the operators
\bes T \left( \suk c_k e_k \right):= \lambda \suk c_{k+1}e_k, \, \,U \left( \suk c_k e_k \right):= \lambda^{-1} \suk c_{k}e_{k+1}, \, \, \ctk \in \ltn.\ens
Clearly, for the particular Hilbert space $\ltz$ an indexing of the orthonormal basis as $\etk$ is not natural. An alternative
would be to use the orthogonal decomposition
\bes \ltz = \ell^2(\dots, -2,-1,0)\oplus \ell^2(1,2,3, \dots)
= \ell^2(\mz_{-}) \oplus \ell^2( \mz_{+}) \ens and then apply the above methods to each
of the spaces $ \ell^2(\mz_{-})$ and $ \ell^2(\mz_{+});$ that would then lead to an approximation of
a given frame for $\ltz$ using a union of two suborbits
of appropriately chosen operators $T_1, T_2.$ We will
not go into details with such constructions here but just notice
that they would follow the line of the approach we take in
the next section for the Hilbert space $\ltr.$}
\end{rem}

\subsection{Approximate frame representations in $\ltr.$} \label{191018c}

The key reason that makes the approach in Section \ref{191018b} work
is that the shift operator $T$ in \eqref{231018d} removes the first coordinate
in a given vector in $\ltn$. Thus, the method can not directly
be generalized to the shift operator  on $\ltr.$ For $\ltr$
we could alternatively apply the procedure outlined in Remark \ref{231018c}, but again the indexing of an orthonormal basis
for $\ltr$ as $\etk$ might not appear natural. In this
section we will outline an approach that is more similar to
what we did in Section \ref{191018b}. In the entire section we will
use the following setup:

\vn{\bf General setup:} Consider a frame $\ftk$ for $\ltr$ consisting of
compactly supported functions and with a common upper
bound on the length of the support, see \eqref{071018b} below. Let $\gtk$ denote the
collections of functions from $\ftk$ that are supported
within $[0, \infty[,$ and denote the remaining functions
by $\htk$; note that it is  a consequence
of the frame condition that the families $\gtk$ and $\htk$ indeed
are infinite.
Choose for $k\in \mn$ the scalars
$a(k), b(k),c(k), d(k)\in \mr$ such that
\bee \label{071018c} \supp \, g_k \subseteq [a(k), b(k)],
\, \, \supp \, h_k \subseteq [c(k), d(k)].\ene
Assume  that
\bee \label{071018b} L:= \sup_{k\in \mn}\left(  b(k)-a(k), d(k)-c(k)\right) < \infty. \ene
Fix some $\lambda>1,$ and define the following truncated and
scaled translation operators on $\ltr$:

\bee \label{071018a} U_1f(x)= \lambda^{-1}f(x-1), \ \ & \ &
T_1f(x)= \lambda f(x+1) \chi_{[0,\infty[}(x) \\
 U_2f(x)= \lambda^{-1}f(x+1), \ \ & \ &
T_2f(x)= \lambda f(x-1) \chi_{]-\infty,L]}(x) \ene
Finally, choose increasing sequences $\{\alpha(k)\}_{k=1 }^\infty, \{\gamma(k)\}_{k=1 }^\infty$ of positive numbers such that $\alpha(k+1)-\alpha(k)\ge b(k)$ and
$\gamma(k+1)- \gamma(k) \ge L-c(k)$ for all $k\in \mn,$
and let
\bee \label{071018d} \varphi_1:= \sum_{n=1}^\infty U_1^{\alpha(n)}g_n, \, \,  \varphi_2:= \sum_{n=1}^\infty U_2^{\gamma(n)}h_n. \ene
The standing assumptions imply that for $n\neq \ell$ we have
$U_1^{\alpha(n)}g_n \bot U_1^{\alpha(\ell)}g_\ell$ and
$U_2^{\gamma(n)}h_n \bot U_2^{\gamma(\ell)}h_\ell;$
thus, the same estimate as in \eqref{031218a} applies
to show that the vectors $\varphi_1$ and $\varphi_2$ are well-defined.
We will now show how the frame $\ftk= \gtk \cup \htk$ can be approximated
using suborbits of the two operators $T_1$ and $T_2,$ generated by the
functions $\varphi_1,$ respectively $\varphi_2:$

\bt \label{071018f} In the above setup, for all $k\in \mn$
it holds that
\bee \label{071018g} || g_k - T_1^{\alpha(k)} \varphi_1||^2
\le \frac{B\lambda^2}{\lambda^2-1} \, \lambda^{-2 \, [\alpha(k+1)-\alpha(k)]}\ene
and
\bee \label{071018h} || h_k - T_2^{\gamma(k)} \varphi_2||^2
\le \frac{B\lambda^2}{\lambda^2-1} \, \lambda^{-2 \, [\gamma(k+1)-\gamma(k)]}.\ene
In particular, if we for some $\epsilon \in ]0, A[$ choose
the sequences  $\{\alpha(k)\}_{k=1 }^\infty, \{\gamma(k)\}_{k=1 }^\infty$ such that $\alpha(1)= \gamma(1)=0$ and
for $k\in \mn,$
\bee \label{071018k}
\alpha(k+1)-\alpha(k) \ge \max\left(b(k), \frac{k \,\ln(2)+\ln \left(\frac{2B}{\epsilon}\right) + \ln \left( \frac{\lambda^2}{\lambda^2-1}\right)  }{2\ln (\lambda)}\right), \ene and
\bee \label{071018j}
\gamma(k+1)-\gamma(k) \ge \max\left(L-c(k), \frac{k \,\ln(2)+\ln \left(\frac{2B}{\epsilon}\right) + \ln \left( \frac{\lambda^2}{\lambda^2-1}\right)  }{2\ln (\lambda)}\right), \ene then $\{T_1^{\alpha(k)}\varphi_1\}_{k=1}^\infty\cup
\{T_2^{\gamma(k)}\varphi_2\}_{k=1}^\infty$ is an
$\epsilon$-approximation of the frame $\gtk \cup \htk= \ftk$ and
Theorem \ref{181217ens} applies.
\et

\bp The proof follows the lines of the proof of Theorem \ref{011018an}, so we only sketch it. We have already argued that the vectors $\varphi_1$
and $\varphi_2$ are well-defined.
Now, fixing $k\in \mn,$ we have $T_1^{\alpha(k)}U_1^{\alpha(k)}g_k=g_k.$ Also, for
$n<k,$ $\alpha(k)-\alpha(n) \ge \alpha(n+1)-\alpha(n)\ge b(n);$
it follows that then $T_1^{\alpha(k)}U_1^{\alpha(n)}g_n=
T_1^{\alpha(k)-\alpha(n)}g_n=0.$ Now exactly the same argument
as in \eqref{031218b} yields the conclusion in
\eqref{071018g}. The proof of \eqref{071018j} is similar. \ep

We will now show explicitly how Theorem \ref{071018f} can be applied to Gabor frames. Recall that
a Gabor frame is a frame for $\ltr$ having the form
$\mts:= \{e^{2\pi i bx}g(x-na)\}_{m,n\in \mz}$ for a fixed choice
of a function $g\in \ltr$ and some $a,b>0.$ In order to apply Theorem \ref{071018f}
we only need to assume that $g$ is compactly supported. Furthermore, for
special choices of the involved parameters we can be completely explicit
about how to choose suitable powers $\alpha(k)$ and $\gamma(k)$ in
Theorem \ref{071018f}:
\bex
Consider a Gabor frame $\mts$  and assume that  $\supp\, g\subseteq [0,C]$ for some $C\in\mr$.
Choose the sequence $\gtk$ to consist of  only the positive  translates $\{E_{mb}T_{na}g\}_{m\in\mz,n\in\mn_0}$. Furthermore, we order $\gtk$ such that
$g_1:=g$ and, in general, if
$g_k=E_{m'b}T_{n'a}g$ for some $m'\in\mz, n'\in\mn_0,$ then $g_{k+1}$  corresponds to one of the following functions:
\bes E_{(m'\pm1)b}T_{n'a}g, \quad E_{m'b}T_{(n'\pm1)a}g.
\ens In other words, the ordering is chosen such that going from any $g_k$ to $g_{k+1},$
the translation parameter changes by at most $a.$ Thus, we have $\supp \, g_{k+1}\subseteq (\supp \, g_k + a)\cup (\supp \,g_k -a)\cup \, \supp \, g_k.$ Combined with
the fact that all the functions $g_k$ are supported within the positive real axis, this
implies that
for each $k\in \mn$ there
exists a  positive integer $\ell_k\le k$ such that $\supp \,g_k\subseteq [(\ell_k-1)a,C+(\ell_k-1)a]$.
% Therefore we can take $b(k)=C+(k-1)a$ in \eqref{071018c}.
In a similar fashion, we let $\htk$
consist of the subsequence  $\{E_{mb}T_{na}g\}_{m\in\mz,n\in \{\dots, -3,-2,-1\}},$
ordered such that $h_1=T_{-a}g$ and  the translation parameters in $h_{k+1}$ and $h_k$ differ with at most $a$. Thus
$\supp \, h_k\subseteq [-r_ka,C-r_ka ]$ for some $r_k\in\mn$ and $r_k\leq k$. %Thus we can take $c(k)=-ka$ in \eqref{071018c}.

Now, consider the case $\lambda=\sqrt{2}$ and choose
the ``tolerance" as  $\epsilon=2^{-j}$
for some $j\in \mn.$ Furthermore, choose an upper frame bound
of the form  $B=2^N$ for a sufficiently large value of $N\in\mn$. Then, let
$\alpha(1)= \gamma(1)=0$ and define the sequences $\{\alpha(k)\}_{k=1}^\infty$ and  $\{\gamma(k)\}_{k=1}^\infty$ recursively by
\bes
\alpha(k+1)-\alpha(k)&=&C+(a+1)k- a + N+j+2,\quad k\in\mn.\\
\gamma(k+1)-\gamma(k)&=&C+(a+1)k+N+j+2,\quad k\in\mn.
\ens
Explicitly,  this gives
\bes
\alpha(k)&=&(k-1)\left[k\frac{(a+1)}{2} + C-a+N+j+2\right],\quad k\in\mn, \\
\gamma(k)&=&  (k-1)\left[k\frac{(a+1)}{2} + C+N+j+2\right],\quad k\in\mn.
\ens

Then, defining $\varphi_1,\varphi_2$ as  in \eqref{071018d}.  Theorem~\ref{071018f} implies that the sequence $\{T^{\alpha(k)}\varphi_1\}_{k=1}^\infty\cup \{T^{\gamma(k)}\varphi_2\}_{k=1}^\infty$ is a frame for $L^2(\mr)$ and it is an $\epsilon$-approximation of  the frame  $\{E_{mb}T_{na}g\}_{m,n\in\mz}$. \ep
\enx

\section*{Appendix: Multi-generator representations}
In contrast to the rest of the paper, the
approximate frame representation in Theorem \ref{071018f} is using
suborbits of two operators rather than just one. There is no detailed account of
frame representations using multiple orbits in the literature;
for this reason we
end the paper with a short collection of results concerning
frame representations using multiple orbits. 

\bt \label{80923fn}   Let $\ftk$ be a frame, which is
norm-bounded below. Then there is a finite collection of
vectors from $\ftk,$ to be called $\varphi_1, \dots, \varphi_J,$ and
corresponding bounded operators $T_j: \h \to \h$, such that
\bee \label{50614a} \ftk = \bigcup_{j=1}^J \{T_j^n \varphi_j\}_{n=0}^\infty.\ene
\et

Note that in contrast with the single-orbit representation in
Theorem \ref{80923f}(i), the representation considered in Theorem \ref{80923fn}
is available for the standard frames considered in applied harmonic analysis, e.g., the frames of translates, Gabor frames, wavelet frames, and generalized shift-invariant systems
generated by a finite collection of nonzero functions.
Note also that while frames with positive and finite excess never
have a representation as an orbit of a bounded operator,
they always have a  representation as a finite union of operator orbit after removal of irrelevant zero-vectors:

\bc Let $\ftk$ be a frame with finite excess
and assume that $f_k \neq 0$ for all $k\in \mn.$
Then $\ftk$ has a representation as in
\eqref{50614a}.\ec

It is natural to ask how large one has to choose
the parameter $J$ in Theorem \ref{80923fn}, and how
$J$ is related to the redundancy and excess of
the frame $\ftk.$
First,
by Theorem \ref{80923f}(i) the condition
of $\ftk$ being norm-bounded below  implies that for
overcomplete frames 
we necessarily have $J\ge 2$ in Theorem \ref{80923fn}.
Technically, Theorem \ref{80923fn} was proved
using the Feichtinger Theorem (an equivalent form of
the long-standing Kadison--Singer Conjecture, which
was finally confirmed in \cite{MSS}) which states that
every frame that is norm-bounded below can be represented
as a finite union of Riesz sequences.
Thus, it would be natural to expect that we would have to choose $J$ at least
as large as the minimal number of Riesz sequences in
the Feichtinger decomposition. This turns out not to be
true. We will demonstrate this by exhibiting a class
of frames for which we can
actually take $J=2,$ regardless of the number of
Riesz sequences appearing in the Feichtinger decomposition.

\bex Consider any frame that has
a representation as a bi-infinite orbit
$\{T^n \varphi\}_{n\in \mz}$ for a bounded
and bijective operator
$T: \h \to \h.$ Then
\bee \label{80923h} \{T^n \varphi\}_{n\in \mz}=
\{T^n \varphi\}_{n=0}^\infty \cup
\{(T^{-1})^n T^{-1}\varphi\}_{n=0}^\infty,\ene
i.e., the frame  $\{T^n \varphi\}_{n\in \mz}$ can be represented as a union of two orbits.
Frames with a representation as a bi-infinite
orbit of a bounded operator were characterized in
\cite{olemmaarzieh-2}; in particular, the assumptions
are satisfied for frames
of translates $\{T_{kb}\varphi\}_{k\in \mz}$ in $\ltr.$
Such
frames can be constructed with arbitrarily large
redundancy (thus forcing a large number $J$ in the
Feichtinger decomposition). Indeed,  it is well known that
the integer-translates $\{T_k \mbox{sinc}\}_{k\in \mz}$ of the sinc-function
forms an ortonormal
basis for the Paley-Wiener space
\bes PW:= \{f\in \ltr \, \big| \, \supp \hat{f} \subseteq [-1/2,1/2]\}.\ens
It follows that for any $N\in \mn,$ the family
$\{T_{k/N} \mbox{sinc}\}_{k\in \mz} $ is an overcomplete
frame for $PW$ consisting of the union of $N$
orthonormal bases. Thus the proof of Theorem \ref{80923f} (ii)
suggests that we must take  $J\ge N;$ however, by
\eqref{80923h}  the frame can be decomposed
as a union of just two operator orbits. \ep \enx

An intermediate step between  a representation of a frame
as a single orbit and as in \eqref{50614a} would be to consider
multiple orbits generated by the same operator $T.$ The next example shows that  in contrast to the case of frames represented as a single orbit, such frames are not forced to be linearly independent, and they might be norm-bounded below.
The example  consists of a construction of two families of vectors, each of which are complete and satisfy the Bessel condition but
not the lower frame condition; however, the union of the two systems forms a frame
that can be represented as a union of two operator systems, generated by the
same operator $T$ but using different generators.

\bex Let $\etk$ denote an orthonormal basis for $\h.$ Our purpose
is to consider the sequence $\{e_k + e_{k+1}\}_{k=1}^\infty \cup \{e_k - e_{k+1}\}_{k=1}^\infty $ and show that
\bei
\item[(i)] $\{e_k + e_{k+1}\}_{k=1}^\infty \cup \{e_k - e_{k+1}\}_{k=1}^\infty $ is
linearly
dependent, and hence can not be represented as an orbit of a single
operator (bounded or not).
\item[(ii)] $\{e_k + e_{k+1}\}_{k=1}^\infty \cup \{e_k - e_{k+1}\}_{k=1}^\infty $ is
a frame and can be represented in the form $\cup_{j=1}^2\{T^n \varphi_j\}_{n=0}^\infty$
for a bounded operator $T: \h \to \h$ and some $\varphi_1, \varphi_2 \in \h;$
\eni

First, the result in (i) follows immediately from the
observation that
$0= (e_1+e_2)- (e_1-e_2) - (e_2+e_3) - (e_2-e_3).$
Hence, the family  $\{e_k + e_{k+1}\}_{k=1}^\infty \cup \{e_k - e_{k+1}\}_{k=1}^\infty $  can not be represented as an orbit of a single
operator. To prove that it can be represented as a union of two
operator orbits we note that
by Example 5.4.6 in
\cite{CB} the
family  $\{e_k + e_{k+1}\}_{k=1}^\infty $ is a complete Bessel sequence, but
not a frame. Considering the bounded linear operator $T:\h \to \h$ defined by
$Te_k:=e_{k+1}, \, k\in \mn$ and letting $\varphi_1:= e_1+e_2,$ we have
$\{e_k + e_{k+1}\}_{k=1}^\infty= \{T^n \varphi_1\}_{n=0}^\infty.$ A slight modification of the
argument proves that also $\{e_k - e_{k+1}\}_{k=1}^\infty$ is a complete Bessel sequence but
not a frame; also, with the same operator $T$ and $\varphi_2:= e_1-e_2,$ we have
$\{e_k - e_{k+1}\}_{k=1}^\infty= \{T^n \varphi_2\}_{n=0}^\infty.$

To prove the frame property in (ii) it is enough to find
a Bessel sequence leading to perfect reconstruction. By a direct calculation,
for any $f\in \h,$
\bes f= \suk \la f, \frac12 e_k\ra (e_k + e_{k+1})+ \suk \la f, \frac12 e_k\ra (e_k - e_{k+1});\ens
this proves that $\{e_k + e_{k+1}\}_{k=1}^\infty \cup \{e_k - e_{k+1}\}_{k=1}^\infty $ indeed is a frame, with dual frame
$\{\frac12 e_k \}_{k=1}^\infty \cup
\{\frac12 e_k \}_{k=1}^\infty$.
\ep \enx

While certain frames with finite excess can be represented
using unions of operator orbits, this possibility disappear
if we insist on all orbits being generated by the same
operator:

\bpr \label{80423a} Consider a frame $\ftk$ with finite and strictly positive excess. Then there does not exist a bounded operator $T: \h \to \h$ and a finite number of
generators $\{\varphi_j\}_{j=1}^J$ such that
\bee\label{180518a}
\ftk = \bigcup_{j=1}^J \{ T^n \varphi_j\}_{n=0}^\infty.
\ene \epr

\bp Let $\ftk$ be a frame for $\h$ with finite excess. Assume that there is an operator $T\in B(\h)$ and $\varphi_1,\dots,\varphi_J\subset\ftk$ such that \eqref{180518a} holds.
Since $\ftk$ has positive excess, there exists an $\ell\in\mn_0$ and $j_0\in \{1,\dots,J\}$ and scalars $\{c_n^j\}_{j=1,\cdots,J,n\in\mn_0}\subset\mc$ such that
\[T^\ell \varphi_{j_0}= \sum_{j\neq j_0} \sum_{n=0}^\infty c_n^j T^n\varphi_j + \sum_{n\neq \ell} c_n^{j_0} T^n\varphi_{j_0}
\] Using that the frame $\ftk$ is assumed to have finite excess, we can choose 
$N\in\mn$ such that $\cup_{j=1}^J \{ T^n \varphi_j\}_{n=N}^\infty$ is a Riesz sequence. Now, take any $M>N$. Then, acting with $T^M$ yields
\bes T^{M+\ell } \varphi_{j_0}&=& \sum_{j\neq j_0} \sum_{n=0}^\infty c_n^j T^{n+M}\varphi_j + \sum_{n\neq \ell} c_n^{j_0} T^{n+M}\varphi_{j_0}\\
&=&\sum_{j\neq j_0} \sum_{n=M}^\infty c_n^j T^{n}\varphi_j + \sum_{n\geq M, n\neq M+\ell} c_n^{j_0} T^{n}\varphi_{j_0}.
\ens
This implies that
\[ T^{M+\ell } \varphi_{j_0}\in \span \left\{ T^n \varphi_j  \right\}_{n\geq M , j\neq j_0 }
\cup \{T^n \varphi_{j_0} \}_{n\geq M , n\neq \ell + M}\]
contradicting that $\cup_{j=1}^J \{ T^n \varphi_j\}_{n=N}^\infty$ is a Riesz basis. \ep

\vn{\bf Acknowledgment:}  The authors would like to thank the anonymous referees
for their helpful comments which improved the presentation of the results.

\begin{tabbing}
text-text-text-text-text-text-text-text-text-text \= text \kill \\
Ole Christensen \> Marzieh Hasannasab \\
Technical University of Denmark \> Technical University of Kaiserslautern  \\
DTU Compute \> Paul-Ehrlich Stra\ss e Geb\"aude 31 \\
Building 303, 2800 Lyngby \> 67663 Kaiserslautern \\
Denmark \> Germany \\
Email: ochr@dtu.dk \> hasannas@mathematik.uni-kl.de
\end{tabbing}


\begin{thebibliography}{99}
\bibitem{A1} Aldroubi, A.,  Cabrelli, C.,  Molter, U., and Tang, S.:
{\it Dynamical sampling.} Appl. Harm. Anal. Appl.  {\bf 42} (2017), 378--401.
%
%
\bibitem{A2} Aldroubi, A., Cabrelli, C., \c{C}akmak, A. F., Molter, U., and  Petrosyan, A.: {\it Iterative actions of normal operators.} J. Funct. Anal. {\bf 272} no. 3 (2017),
1121--1146.

\bibitem{A3} Aldroubi, A. and Petrosyan, A.: {\it Dynamical sampling and systems from iterative actions of operators}, Frames and other bases in abstract and function spaces, 15–26, Appl. Numer. Harmon. Anal., Birkhäuser/Springer,
     Cham, 2017.




\bibitem{BoSp2} Bownik, M., and Speegle, D.: {\it    Linear independence of Parseval wavelets.} Illinois J. Math. {\bf 54} (2010), 771--785.


\bibitem{CMPP} Cabrelli, C., Molter, U, Paternostro, V., and Philipp, F.: {\it Dynamical
	sampling on finite index sets.} Preprint, 2017.

\bibitem{Chen} Chen, K. Y.: {\it On aperiodicity and hypercyclic weighted translation operators.} J. Math. Anal. Appl. {\bf 462} (2018),  no. 2, 1669--1678


\bibitem{CB} Christensen, O.:
{\it An introduction to frames and Riesz bases.} Second expanded
edition. Birkh\"auser (2016)




\bibitem{olemmaarzieh}
Christensen, O. and Hasannasab, M.:
{\it Frame properties of systems arising via iterative actions of operators}
Appl. Comp. Harm. Anal. {\bf 46} (2019), 664--673.


\bibitem{olemmaarzieh-2}
Christensen, O., and Hasannasab, M.:
 {\it Operator representations of frames: boundedness, duality, and stability}
\newblock Integral Equations an Operator Theory {\bf 88} (2017), no. 4, 483--499.








\bibitem{CHP} Christensen, O., Hasannasab, M., and Philipp, F.: {\it Frame Properties of Operator Orbits.} Accepted for publication in Math. Nach., 2019

 \bibitem{olemmaarzieh-3}
Christensen, O., Hasannasab, M., and Rashidi, E.:
{\it Dynamical sampling and frame representations with bounded
operators.}
\newblock  J. Math. Anal. Appl. {\bf 463} no. 2 (2018), 634--644.







\bibitem{Grivaux}
Grivaux, S.: {\it Construction of operators with prescribed behaviour}, Arch. Math. (Basel), 81 (2003), no. 3, 291--299.

\bibitem{erdmannsurvey}
Grosse-Erdmann K. G.: {\it Universal families and hypercyclic operators} Bulletin of the American Mathematical Society. 1999;36(3):345--81.


\bibitem{erdmann}
Grosse-Erdmann, K. G., and  Manguillot, A. P.: {\it Linear chaos}, Springer Science $\&$ Business Media, (2011).




\bibitem{G5} Gr\"{o}chenig, K.: {\it Localization of frames,
Banach frames, and the invertibility of the frame operator.}
 J. Fourier Anal. Appl. {\bf 10} (2004), 105--132.

\bibitem{kitai}
Halperin, I., Kitai, C., and Rosenthal, P.: {\em On orbits of linear operators.} Journal of the London Mathematical Society 2, {\bf 3} (1985), 561--565.

\bibitem{Linn} Linnell, P.: {\it Von Neumann algebras and linear independence
of translates.} Proc. Amer. Math. Soc. {\bf 127} no. 11 (1999),
3269--3277.


\bibitem{MSS}  Marcus, A. W., Spielman, D. A., and Srivastava, N.: {\it Interlacing families II: Mixed characteristic polynomials and the Kadison--Singer problem.}  Ann. of Math. (2) {\bf 182} no. 1, (2015), 327--350.

\bibitem{rolewicz}
Rolewicz, S.: {\it On orbits of elements.} Studia Mathematica 32, {\bf 1}  (1969), 17--22.






\end{thebibliography}
\end{document}